**The Diophantine equation** $xy = z^n$ **; for** $n = 2, 3, 4, 5, 6.$;

**the Diophantine equation** $xyz = w^2$; **and the**

**Diophantine system** $\begin{Bmatrix} xy = v^2 \\ yz = w^2 \end{Bmatrix}$

Konstantine Zelator

University of Wisconsin Colleges

UW – Marinette

Mathematics

750 W. Bay Shore St.

Marinette, WI 54143

U.S.A.

Also:   Konstantine Zelator

P.O. Box 4280

Pittsburgh, PA 15203

U.S.A.

Email addresses: 1) konstantine_zelator@yahoo.com

## 1. Introduction

In this work, we accomplish three goals. First we determine the entire family of positive integer solutions to the three-variable Diophantine equation,

$x \cdot y = z^n$ ; for $n = 2, 3, 4, 5$ and $6$.

This is done in **Theorems 2** through **6** respectively.

The family of solutions, is in each case described parametrically. Specifically, in the case of $n = 2$, the set of solutions is described in terms of parametric formulas involving three (independent) positive integer parameters.

For $n = 3$, we obtain five-parameter family of solutions.
For $n = 4$, also a 5-parameter family of solutions.
For $n = 5$, a 7-parameter family of solutions.
Likewise for $n = 6$, also a seven-parameter family of solutions.

The second goal of this paper, is determining the solution set, in positive integers; of the four-variable Diophantine equation $xyz = w^2$.

This is done in **Theorem 7**.

The solution set is described in terms of formulas involving six (independent) positive integer parameters.

Finally, the third goal is solving the Diophantine,

5-variable system; $\begin{cases} xy = v^2 \\ yz = w^2 \end{cases}$

We determine all the positive integer solutions of this system. This is done in **Theorem 8**.

The entire solution set is parametrically described in terms of eight parameters.

In this last section of this work (section 12), we mention the obvious geometric interpretations. Each of the above equations (as well as the 2-equation system) can be interpreted geometrically in an obvious way. For example, think of the equation $xyz = w^2$. One can try to determine some or all 3-D boxes; that have integer side lengths and volume which is an integer or perfect square.

After searching the internet, this author has not been able to find material on the equations of this paper. However, W. Sierpinski's book, *Elementary Theory of Numbers* (see reference **[1]**), list two papers with material on equations of the form, $x_1 \cdot x_2 \cdots x_n = t^k$. One is a 3-page paper published in 1955 (see reference **[3]**). And the other, a 10-page paper published in 1933 (see reference **[4]**).

This author has not been able to access these two papers. Even if some of the results in those two works; overlap with some of the results in this work, it is quite likely, that the methods used in this paper; will be different from those used in the above mentioned papers.

2. **Three Lemmas and Proposition 1**

The first lemma, very well-known and widely used in number theory, is known as Euclid's. From a proof see reference **[1]** or reference **[2]**.

**Lemma 1** (Euclid's Lemma)

*Let a,b,c be positive integers such that the integer is a divisor of the product bc and with a and b being relatively prime. Then a must be a divisor of c.*

The following lemma can be proved by using the Fundamental Theorem of Arithmetic. It can also be proved without the use of the Fundamental Theorem (for example, refer to reference **[1]**)

**Lemma 2**

*Let a, b, and n be positive integers. If $a^n$ is a divisor of $b^n$, then a must be a divisor of b.*

One may use the Fundamental Theorem (factorization of a natural number into a product of prime powers), to establish the lemma below, lemma 3. In reference **[1]**, the reader can find a proof of lemma 3 that makes use of lemmas 1 and 2; but not the factorization theorem of a positive integer into prime powers.

**Lemma 3**

*Suppose that a, b, c, n are positive integers such that $ab = c^n$, and that a and b are relatively prime. Then, $a = a_1^n$, $b = b_1^n$, and $c = a_1 b_1$; for some positive, relatively prime positive integers $a_1$ and $b_1$; and with $a_1 b_1 = c$*

**Proposition 1**

Let k be a positive integer. And a,b,c, and n, positive integers such that $ab = k \cdot c^n$; and with a and b being relatively prime. Then $a = d a_1^n$, $b = K \cdot b_1^n$, and $c = a_1 b_1$; where $a_1$ and $b_1$ are relatively prime positive integers. And $d$ and $K$ are also relatively prime integers such that $d \cdot K = k$

**Proof**

We have $ab = kc^n$     **(i)**

Let $d = \gcd(a, k)$; the greatest common divisor of $a$ and $k$.

Then  $a = A \cdot d$ and $k = K \cdot d$,

Where $A$ and $K$ are relatively prime positive

integers; $\gcd(A, K) = 1$     **(ii)**

From **(i)** and **(ii)** we obtain

$$A \cdot b = K \cdot c^n \quad \textbf{(iii)}$$

Since $\gcd(A, K) = 1$. It follows by Lemma 1 and **(iii)** that $K$ must be a divisor $b$; and so,

$b = K \cdot B$     **(iv)**

for some positive integer $B$.

Combining **(iv)** with **(iii)** yields,

$$A \cdot B = c^n \quad \textbf{(v)}$$

Since $A$ is a divisor of $a$; and $B$ is a divisor of $b$; and $\gcd(a,b) = 1$. It follows that $A$ and $B$ are relatively prime. Thus equation **(v)** and Lemma 3 imply,

$$A = a_1^n \text{ and } B = b_1^n \quad \textbf{(vi)}$$

where $a_1$ and $b_1$ are relatively prime positive integers. And so, $a_1^n \cdot b_1^n = c^n$; which gives $c = a_1 b_1$.

Further we have $\gcd(K, d) = 1$. This follows from the fact that $d$ is a divisor of $A$; $K$ is a divisor of $b$; and the hypothesis $\gcd(a,b) = 1$. Finally from **(iv)**, **(vi)** and **(ii)**, we have $a = d \cdot a_1^n$ and $b = K \cdot b_1^n$. ∎

3. **Theorem 1 and its proof**

**Theorem 1**

*Let k be a fixed (or given) positive integer, and consider the 3-variable Diophantine equation,*

$xy = k \cdot z^n$ *(n a positive integer)*

*Then all the positive integer solutions can be described by the parametric formulas (two parameters $t_1$ and $t_2$) $x = k_1 \cdot t_1^n$, $y = k_2 \cdot t_2^n$, $z = t_1 t_2$.*

*Where $t_1$ and $t_2$ can be any two relatively prime positive integers. And $k_1$, $k_2$ are relatively prime positive integers such that $k_1 k_2 = k$.*

**Proof**

If the positive integers $x, y, z$ satisfy the parametric formulas stated in **Theorem 1**. Then a straightforward calculation shows that the ordered triple $(x, y, z)$ is a positive integer solution to the given equation. Then by **Proposition 1** it follows that $x$, $y$ and $z$ must satisfy the above parametric formulas (with $d = k_1$ and $K = k_2$; $a_1 = t_1$ and $b_1 = t_2$). ∎

4. **Proposition 2 and its Proof**

**Proposition 2**

*Let n be a positive integer, n ≥ 2. Consider the three-variable Diophantine equation,*
$xy = z^n$ *(x,y,z are positive integer variables)*
*Then, this equation is equivalent to the six-variable Diophantine system*

$$\begin{cases} d^2 \cdot X \cdot Y = Z^n \cdot w^{n-2} \\ w^{n-2} = vd^2 \end{cases}$$

*Where $d, X, Y, Z, w, v$ are positive integer variables such that $x = \delta \cdot X, y = \delta \cdot Y, z = w \cdot Z, \delta = w \cdot d$*
*And $\delta = \gcd(x, y), \gcd(X, Y) = 1, w = \gcd(z, \delta), \gcd(Z, d) = 1$*

**Proof**

We start with, $xy = z^n$, $n \geq 2$.     **(1)**
And let $\delta = \gcd(x, y)$. So that,

$$\begin{cases} x = \delta \cdot X, y = \delta \cdot Y \\ \gcd(X, Y) = 1 \end{cases} \quad (2)$$

X and Y positive integers

By **(1)** and **(2)** we obtain,

$$\delta^2 \cdot X \cdot Y = z^n \quad (3)$$

Let $w = \gcd(\delta, z)$. So that,

$$\begin{cases} \delta = w \cdot d \text{ and } z = w \cdot Z, \\ \text{for positive integers } d \text{ and } Z; \\ \text{and with } \gcd(d, Z) = 1 \end{cases} \quad (4)$$

From **(3)** and **(4)** we further get,

$d^2 \cdot w^2 \cdot X \cdot Y = w^n \cdot Z^n$;
And since $n \geq 2$; $d^2 \cdot X \cdot Y = w^{n-2} \cdot Z^n$

The proof is complete ■

5. **Theorems 2 and 3; and their proofs**

We state **Theorem 2.**

**Theorem 2**

*Consider the 3-variable Diophantine equation $xy = z^2$. All the positive integer solutions $(x, y, z)$ are given by the three parameter (parameters $\delta$, $p_1$ and $p_2$) formulas,*

$$x = \delta \cdot p_1, \quad y = \delta \cdot p_2, \quad z = \delta p_1 p_2$$

*Where $\delta, P_1, P^2$ can be any positive integers such that $\gcd(P_1, P_2) = 1$*

**Proof:** If $x, y, z$ satisfy the stated formulas. A straight forward calculation shows that $(x, y, z)$ is a solution. Next, the converse. Assume $(x, y, z)$ to be a solution. We apply Proposition 2 with $n = 2$. The given equation is equivalent to the Diophantine system

$$\begin{cases} d^2 \cdot X \cdot Y = Z^2 \\ 1 = v \cdot d^2 \end{cases} \quad (5)$$

Where $d, X, Y, Z, w, v$ are positive integer variables such that $x = \delta \cdot X$, $y = \delta \cdot Y$, $z = w \cdot Z$, $\delta = w \cdot d$, $\delta = \gcd(x, y)$, $\gcd(X, Y) = 1$, $w = \gcd(z, \delta)$, $\gcd(Z, d) = 1$

From **(5)** we get $v = 1 = d$ and,

$$X \cdot Y = Z^2 \quad (6)$$

Since $X$ and $Y$ are relatively prime. Equation (6) and Lemma 3 imply,

$$\begin{cases} X = p_1^2 \text{ and } Y = p_2^2; \text{ and so } Z = p_1 p_2 \\ \text{where } p_1, p_2 \text{ are relatively prime positive integers; and so } Z = p_1 p_2 \end{cases} \quad (7)$$

From (7), $d = 1$, $w = \delta$, $z = w \cdot Z$, $x = \delta \cdot X$, $y = \delta \cdot Y$.

We obtain $x = \delta \cdot p_1^2$, $y = \delta \cdot p_2^2$, and $z = \delta \cdot p_1 p_2$ ∎

**Theorem 3**

*Consider the 3-variable Diophantine equation,*

$$xy = z^3.$$

*All the positive integer solutions of this equation can be described by the five-parameter formulas (parameters $d, v_1, v_2, m, l$), $x = d^3 \cdot v_1^2 \cdot v_2 \cdot m^3$, $y = d^3 \cdot v_1 \cdot v_2^2 \cdot l^3$, $z = d^2 \cdot v_1 \cdot v_2 \cdot m \cdot l$ Where $d, v_1, v_2, m, l$ can be any positive integers such that $\gcd(v_1, v_2) = 1 = \gcd(m, l)$*

**Proof**

If $x, y, z$ satisfy the stated formulas. Then a straightforward calculation verifies that $(x, y, z)$ is a solution. Next, the converse.

By Proposition 2 with $n=3$, the given equation is equivalent to the system,

$$\begin{cases} d^2 \cdot X \cdot Y = Z^3 \cdot w \\ w = v \cdot d^2 \\ \text{Where } d, X, Y, Z, w, v \text{ are positive integer variables satisfying} \\ x = \delta \cdot X,\ y = \delta \cdot Y,\ z = w \cdot Z,\ \delta = w \cdot d;\ \text{and,} \\ \delta = \gcd(x,y),\ \gcd(X,Y) = 1,\ w = \gcd(z,\delta),\ \gcd(Z,d) = 1 \end{cases} \quad (8)$$

From **(8)** we further obtain,

$$X \cdot Y = v \cdot Z^3 \quad (9)$$

Since $X$ and $Y$ are relatively prime, equation (9) and Proposition 1 imply that,

$$\begin{cases} X = v_1 \cdot m^3,\ Y = v_2 \cdot l^3,\ Z = m \cdot l \\ \text{where } v_1, v_2, m, l \text{ are positive integers such that,} \\ \gcd(v_1, v_2) = 1 = \gcd(l, m) \text{ and } v_1 v_2 = v \end{cases} \quad (10)$$

Accordingly from **(10)** and **(8)** we get,

$$w = v_1 v_2 d^2 \text{ and } \delta = w \cdot d = v_1 v_2 d^3 \quad (11)$$

Hence from **(10)**, **(8)**, and **(11)** we further have,

$x = d^3 \cdot v_1^2 \cdot v_2 \cdot m^3$, $y = d^3 \cdot v_1 \cdot v_2^2 \cdot l^3$, and $z = d^2 \cdot v_1 \cdot v_2 \cdot m \cdot l$.

The proof is complete. ∎

### 6. Theorem 4 and its Proof

**Theorem 4**

*Consider the 3-variable diophantine equation,*
$$xy = z^4$$

*All the positive integer solutions of this equation can be described by the five-parameter formulas (parameters $d, t_1, t_2, f, g$), $x = d^2 \cdot t_1^3 \cdot t_2 \cdot f^4$, $y = d^2 \cdot t_1 \cdot t_2^3 \cdot g^4$, $z = d \cdot t_1 \cdot t_2 \cdot f \cdot g$*

*Where $d, t_1, t_2, f, g$ can be any positive integers such that $\gcd(t_1, t_2) = 1 = \gcd(f, g)$*

**Proof**

If $x, y, z$ satisfy the stated parametric formulas. Then an easy calculation verifies that the triple $(x, y, z)$ is a solution. Now, suppose that $(x, y, z)$ is a positive integer solution. We will show that it must have the form described by the above formulas.

By Proposition 2 ; equation $xy = z^4$ is equivalent to the system (the case $n = 4$)

$$\begin{cases} d^2 \cdot X \cdot Y = Z^4 \cdot w^2 \\ w^2 = v \cdot d^2 \\ \text{Where } d, X, Y, Z, w, v \text{ are positive integer variables such that,} \\ x = \delta \cdot X, \ y = \delta \cdot Y, \ z = w \cdot Z, \ \delta = w \cdot d, \\ \delta = \gcd(X,Y), \ \gcd(X,Y) = 1, \ w = \gcd(z,d), \ \gcd(Z,d) = 1 \end{cases} \quad (12)$$

From $w^2 = v \cdot d^2$ in **(12)**; $d^2$ is a divisor of $w^2$. Thus, by lemma 2; $d$ must be a divisor of $w$. We have
$$w = t \cdot d, \text{ for some positive integer } t \quad (13)$$

Combining **(13)** and **(12)** yields,
$$\begin{cases} X \cdot Y = Z^4 \cdot v \\ \text{and } v = t^2 \end{cases} \text{; or equivalently,}$$

$$XY = Z^4 \cdot t^2 = (Z^2 \cdot t)^2 \quad (14)$$

Since $X$ and $Y$ are relatively prime. Lemma (3) and equation (14) imply that,
$$\begin{cases} X = F^2, \ Y = G^2, \text{ and } Z^2 \cdot t = F \cdot G \\ \text{Where } F \text{ and } G \text{ are positive integers} \\ \text{which are relatively prime; } \gcd(F,G) = 1 \end{cases} \quad (15)$$

Since $\gcd(F,G) = 1$. The third equation in **(15)** together with Proposition 1 imply that,
$$\begin{cases} F = t_1 \cdot f^2, \ G = t_2 \cdot g^2, \ Z = f \cdot g; \\ \text{where } t_1, t_2, f, g \text{ are positive integers such that} \\ \gcd(t_1, t_2) = 1 = \gcd(f,g) \text{ and } t_1 t_2 = t \end{cases} \quad (16)$$

From $t_1 t_2 = t$, $v = t^2 = t_1^2 t_2^2$ (see above), $w = t \cdot d = t_1 \cdot t_2 \cdot d$.

We further get [go back to **(12)**] $\delta = w \cdot d = t_1 t_2 d^2$. And so,

$x = \delta \cdot X = t_1 t_2 d^2 \cdot F^2 = t_1 t_2 d^2 \cdot (t_1 f^2)^2 = d^2 \cdot t_1^3 \cdot t_2 \cdot f^4$

$y = \delta \cdot X = t_1 t_2 d^2 \cdot G^2 = t_1 t_2 d^2 (t_2 \cdot g^2)^2 = d^2 \cdot t_1 \cdot t_2^3 \cdot g^4$

$z = w \cdot Z = t_1 t_2 d \cdot f \cdot g = d \cdot t_1 \cdot t_2 \cdot f \cdot g$ ∎

## 7. Theorem 5 and its proof

*Theorem 5*

*Consider the 3-variable Diophantine equation,*
$$xy = z^5$$
*All the positive integer solutions of this equation can be described by then seven-parameter formulas*

(parameters $q, e_1, e_2, r_1, r_2, i_1, i_2$),

$x = q^5 \cdot e_1^3 \cdot e_2^2 \cdot r_1^4 \cdot r_2 \cdot i_1^5$, $y = q^5 \cdot e_1^2 \cdot e_2^3 \cdot r_1 \cdot r_2^4 \cdot i_2^5$,

and $z = q^2 \cdot e_1 \cdot e_2 \cdot r_1 \cdot r_2 \cdot i_1 \cdot i_2$

Where $q, e_1, e_2, r_1, r_2, i_1, i_2$ can be any positive integers

such that $\gcd(e_1, e_2) = \gcd(r_1, r_2) = \gcd(i_2, i_2) = 1$

**Proof**

If $x, y, z$ satisfy the stated parametric formulas. Then an easy calculation verifies that $(x, y, z)$ is a solution. Now the converse, suppose that $(x, y, z)$ is a solution. By Proposition 1, with $n = 5$. The above equation is equivalent to the Diophantine system,

$$\begin{cases} d^2 \cdot X \cdot Y = Z^5 \cdot w^3 \\ w^3 = v \cdot d^2 \\ \text{Where the positive integer variables } d, X, Y, Z, w, v \text{ satisfy,} \\ x = \delta \cdot X, \ y = \delta \cdot Y, \ z = w \cdot Z, \ \delta = w \cdot d; \text{ and} \\ \delta = \gcd(x, y), \ \gcd(X, Y) = 1, \ w = \gcd(z, \delta), \ \gcd(Z, d) = 1 \end{cases} \quad (17)$$

From the first two equations in **(17)**, we obtain

$$X \cdot Y = v \cdot Z^5 \quad (18)$$

Since X and Y are relatively prime; equation **(18)** and **lemma 3** imply that,

$$\begin{cases} X = v_1 \cdot i_1^5, \ Y = v_2 \cdot i_2^5, \ Z = i_1 i_2 \\ \text{Where } v_1, v_2, i_1, i_2 \text{ are postitive integers} \\ \text{such that } \gcd(v_1, v_2) = 1 = \gcd(i_1, i_2) \\ \text{And } v_1 v_2 = v \end{cases} \quad (19)$$

We go back to the second equation in **(17)**.
Let $D = \gcd(w, d)$. Then,

$$\begin{cases} w = D \cdot r, \ d = D \cdot q \\ \text{Where } r \text{ and } q \text{ are relatively} \\ \text{prime positive integers; } \gcd(r, q) = 1 \end{cases} \quad (20)$$

So, from **(20)** and the second equation in **(17)** we get,

$$D^3 \cdot r^3 = v \cdot D^2 \cdot q^2; \quad (21)$$
$$D \cdot r^3 = v \cdot q^2$$

Since $\gcd(r,q) = 1$, It follows that $\gcd(r^3, q^2) = 1$. Thus, equation (21) and **lemma 1** imply that $q^2$ must be a divisor of D; and so,

$$\begin{cases} D = q^2 \cdot e, \\ \text{for some positive integer } e \end{cases} \quad (22)$$

Combining (22) and (21) gives,

$$v = e \cdot r^3$$

But $v = v_1 v_2$ [see (19)]. And so we further get

$$v_1 v_2 = e \cdot r^3 \quad (23)$$

Since $\gcd(v_1, v_2) = 1$. Proposition and equation (23) imply,

$$\begin{cases} v_1 = e_1 \cdot r_1^3, \\ \text{Where } e_1, e_2, r_1, r_2 \text{ are positive} \\ \text{integers such that,} \\ \gcd(e_1, e_2) = 1 = \gcd(r_1, r_2), \\ e_1 e_2 = e; \text{ and } r_1 r_2 = r \end{cases} \quad (24)$$

We make our way back. From (22) and (24) we have $D = q^2 \cdot e_1 \cdot e_2$.

And so by (20), $d = D \cdot q = q^3 \cdot e_1 \cdot e_2$

Also from (20) and (22), $w = D \cdot r = q^2 \cdot e_1 \cdot e_2 \cdot r_1 \cdot r_2$.

And $\delta = w \cdot d = (q^2 e_1 e_2 r_1 r_2) \cdot (q^3 e_1 e_2) = q^5 \cdot e_1^2 e_2^2 \cdot r_1 r_2$

Further, from (19), $X = v_1 \cdot i_1^5 = e_1 \cdot r_1^3 \cdot i_1^5$ and $Y = v_2 \cdot i_2^5 = e_2 \cdot r_2^3 \cdot i_2^5$

Thus, $x = \delta \cdot X = (q^5 \cdot e_1^2 e_2^2 r_1 r_2)(e_1 r_1^3 i_1^5) = q^5 \cdot e_1^3 e_2^2 \cdot r_1^4 \cdot r_2 \cdot i_1^5$

$x = \delta \cdot X = (q^5 \cdot e_1^2 e_2^2 r_1 r_2)(e_2 r_2^3 i_2^5) = q^5 e_1^2 e_2^3 \cdot r_1 \cdot r_2^4 \cdot i_1^5$

And $z = w \cdot Z = (q^2 e_1 e_2 r_1 r_2)(i_1 i_2) = q^2 e_1 e_2 r_1 r_2 i_1 i_2$ ■

## 8. Theorem 6 and its proof

*Consider the 3-variable Diophantine equation,*

$$xy = z^6$$

*All the positive integer solutions of this equation can be described by the seven-parameter formulas (parameters $p, e_1, e_2, n_1, n_2, j_1, j_2$),*

$$x = p^3 \cdot e_1^4 \cdot e_2^2 \cdot n_1^5 \cdot n_2 \cdot j_1^6, \ y = p^3 \cdot e_1^2 \cdot e_2^4 \cdot n_1 \cdot n_2^5 \cdot j_2^6,$$

and $z = p \cdot e_1 \cdot e_2 \cdot n_1 \cdot n_2 \cdot j_1 \cdot j_2$

Where $p, e_1, e_2, n_1, n_2, j_1, j_2$ can be any positive integers such that
$\gcd(e_1, e_2) = \gcd(n_1, n_2) = \gcd(j_1, j_2) = 1$

**Proof**

If x, y, z satisfy the stated parametric formulas, then a direct calculation verifies that (x, y, z) is a solution. Now the converse, suppose that (x, y, z) is a solution. By **Proposition 1**, with $n=6$. The above equation is equivalent to the Diophantine system,

$$\begin{cases} d^2 \cdot X \cdot Y = Z^6 \cdot w^4 \\ w^4 = v \cdot d^2 \\ \text{Where the positive integer variables } d, X, Y, Z, w, v \text{ satisfy,} \\ x = \delta \cdot X, \ y = \delta \cdot Y, \ z = w \cdot Z, \ \delta = w \cdot d; \text{ and} \\ \delta = \gcd(x, y), \ \gcd(X, Y) = 1, \ w = \gcd(z, \delta), \ \gcd(Z, d) = 1 \end{cases} \quad (25)$$

From **(25)** $\Rightarrow X \cdot Y = v \cdot Z^6$

Also from the second equation in **(25)** ; since $d^2$ is a divisor of $w^4 = (w^2)^2$. We have,

$$\begin{cases} X \cdot Y = v \cdot Z^6 \\ \text{And } w^2 = t \cdot d, \\ \text{for some positive integer } t \end{cases} \quad (26)$$

From the second equation in **(26)** and the second equation in **(25)** we obtain,
$v = t^2$ **(27)**

The first equation in **(26)** implies, by Proposition 1 and since $\gcd(X, Y) = 1$. That,

$$\begin{cases} X = M_1 \cdot j_1^6, \ Y = M_2 \cdot J_2^6, \ Z = j_1 \cdot j_2. \\ \text{Where } M_1, M_2, J_1, J_2; \text{ are} \\ \text{positive integers such that,} \\ \gcd(M_1, M_2) = 1 = \gcd(j_1, j_2) \text{ and} \\ M_1 M_2 = v \end{cases} \quad (28)$$

By **(27)** and $M_1 M_2 = v$, $\gcd(M_1, M_2) = 1$.

We obtain, on account of Lemma 3. That,

$$\left\{ \begin{array}{l} M_1 = m_1^2, \ M_2 = m_2^2, \ t = m_1 m_2 \\ \text{Where } m_1 \text{ and } m_2 \text{ are positive integers such that} \\ \gcd(m_1 m_2) = 1 \end{array} \right\} \quad (29)$$

We return to equation $w^2 = t \cdot d$ in **(26)**
Let $D = \gcd(d, w)$. Then,

$$\left\{ \begin{array}{l} d = p \cdot D, \ w = D \cdot N \\ \text{where p and N are positive integers such that} \\ \gcd(p, N) = 1 \end{array} \right\} \quad (30)$$

So that $w^2 = t \cdot d$ and **(30)** yield,

$$D^2 \cdot N^2 = t \cdot p \cdot D;$$
$$D \cdot N^2 = t \cdot p \quad (31)$$

Since $\gcd(p, N) = 1$; we have $\gcd(p, N^2) = 1$. And so **(31)** implies that p must be a divisor of D:

$$\left\{ \begin{array}{l} D = p \cdot e, \text{ for some} \\ \text{positive integer e} \end{array} \right\} \quad (32)$$

From **(31), (32),** and $t = m_1 m_2$ (in **(29)**),

$$e \cdot N^2 = m_1 m_2 \quad (33)$$

Since $\gcd(m_1 m_2) = 1$. Equation **(33)** and **Proposition 1** imply that,

$$\left\{ \begin{array}{l} m_1 = e_1 \cdot n_1^2, \ m_2 = e_2 \cdot n_2^2, \ e = e_1 e_2, \\ \text{and } N = n_1 n_2 \\ \text{Where } e_1, e_2, n_1, n_2 \text{ are postitive integers} \\ \text{with } \gcd(e_1, e_2) = 1 = \gcd(n_1, n_2) \end{array} \right\} \quad (34)$$

We make our way back using **(34), (32), (30), (29), (28)** and **(25).** We have,

$t = m_1 m_2 = e_1 e_2 n_1^2 n_2^2$, $D = p \cdot e = p e_1 e_2$, $M_1 = m_1^2 = e_1^2 n_1^4$

$w = D \cdot N = p e_1 e_2 n_1 n_2$, $d = p \cdot D = p^2 e_1 e_2$, $M_2 = m_2^2 = e_2^2 n_2^4$

And $\delta = w \cdot d = (p e_1 e_2 n_1 n_2)(p^2 e_1 e_2) = p^3 e_1^2 e_2^2 n_1 n_2$

Next, $X = M_1 \cdot j_1^6 = e_1^2 n_1^4 j_1^6$, $Y = M_2 j_2^6 = e_2^2 n_2^4 \cdot j_2^6$

Thus, $x = X \cdot \delta = (e_1^2 n_1^4 j_1^6)(p^3 e_1^2 e_2^2 n_1 n_2)$;

$x = p^3 \cdot e_1^4 \cdot e_2^2 \cdot n_1^5 \cdot n_2 \cdot j_1^6$, $y = p^3 \cdot e_1^2 \cdot e_2^4 \cdot n_1 \cdot n_2^5 \cdot j_2^6$,

and $z = w \cdot Z = p \cdot e_1 \cdot e_2 \cdot n_1 \cdot n_2 \cdot j_1 \cdot j_2$

∎

### 9. Proposition 3 and its proof

*Let n be a positive integer, $n \geq 2$. Consider the 4-variable Diophantine equation,*

$\quad\quad\quad xyz = w^n \quad\quad$ (x, y, z, w are positive-integer variables)

*Then this equation is equivalent to the seven-variable Diophantine system,*

$$\begin{cases} d^{n-2} = D \cdot v^2 \\ X \cdot Y \cdot z = W^n \cdot D \end{cases}$$

*Where $d, D, v, X, Y, W$ are positive integer variables such that,*

$d = \gcd(w, \delta)$, $w = W \cdot d$, $\delta = d \cdot v$, $\gcd(W, v) = 1$,

$x = \delta \cdot X$, $y = \delta \cdot Y$, $\delta = \gcd(x, y)$, $\gcd(X, Y) = 1$

**Proof**

We start with $xyz = w^n$, $n \geq 2$ \quad\quad\quad (35)

Let $\delta = \gcd(x, y)$. So that,

$\{x = \delta \cdot X, \; y = \delta \cdot Y, \; \gcd(X, Y) = 1\}$ \quad\quad (36)

From (36) and (35) we obtain,

$\delta^2 \cdot X \cdot Y \cdot z = w^n$ \quad\quad\quad (37)

Let $d = \gcd(w, \delta)$. Then,

$\{w = W \cdot d, \; \delta = d \cdot v, \; \gcd(W, v) = 1\}$ \quad\quad (38)

From (38) and (37) we further have

$\delta^2 \cdot X \cdot Y \cdot z = W^n \cdot d^n$;

$d^2 \cdot v^2 \cdot X \cdot Y \cdot z = W^n \cdot d^n$ ;

And since $n \geq 2$; $v^2 \cdot X \cdot Y \cdot z = W^n \cdot d^{n-2}$  (39)

Since $\gcd(v, W) = 1$. It follows that $\gcd(v^2, W^n) = 1$; which together with (39) and Lemma 1; imply that $v^2$ must be a divisor of $d^{n-2}$:

$$\begin{cases} d^{n-2} = D \cdot v^2, \text{ for some} \\ \text{positive integer D} \end{cases} \quad (40)$$

Equations (39) and (40) taken together establish the result. ∎

### 10. Theorem 7 and it proof

**Theorem 7**

*Consider the four-variable Diophantine equation,*

$xyz = w^2$

*All the positive integer solutions of this equation can be described by the six parameter formulas, (parameters $d, r_1, r_2, t, u_1, u_2$)*

$x = d \cdot r_1^2 \cdot u_1$, $y = d \cdot r_2^2 \cdot u_2$, $z = t^2 \cdot u_1 \cdot u_2$,

*and* $w = d \cdot r_1 \cdot r_2 \cdot t \cdot u_1 \cdot u_2$

*And with* $\gcd(r_1, r_2) = 1 = \gcd(u_1, u_2) = 1 = \gcd(t, r_1 r_2)$

**Proof**

If x,y,z,w satisfy the stated parametric formulas, then a straightforward calculation verifies that the quadruple (x,y,z,w) is a solution.
Now the converse, suppose that (x,y,z,w) is a solution. Then, according to **Proposition 3** with $n=2$. We must have,

$$\begin{cases} 1 = D \cdot v^2 \\ X \cdot Y \cdot z = W^2 \cdot D \\ \text{With the positive integer variables } D, v, X, Y, W \\ \text{satisfying,} \\ d = \gcd(w, d), \ w = W \cdot d, \ \delta = d \cdot v, \ \gcd(W, v) = 1, \\ x = \delta \cdot X, \ y = \delta \cdot Y, \ \delta = \gcd(x, y), \ \gcd(X, Y) = 1 \end{cases} \quad (41)$$

From the first equation in (41) we get $D = v = 1$.
And from the second equation $X \cdot Y \cdot z = W^2$

Altogether,
$$\begin{cases} X \cdot Y \cdot z = W^2 \\ w = W \cdot d,\ x = d \cdot X,\ y = d \cdot Y, \\ d = \gcd(X,Y),\ \gcd(X,Y) = 1 \end{cases} \quad (42)$$

Let $\gcd(z, W) = p$. Then,
$$\begin{cases} z = t \cdot p,\ W = r \cdot p, \\ \text{With } \gcd(t, r) = 1 \end{cases} \quad (43)$$

From **(42)** and **(43)** we obtain,

$$X \cdot Y \cdot t = r^2 \cdot p \quad (44)$$

Since $\gcd(t, r) = 1$; it follows that $\gcd(t, r^2) = 1$
And so by **(44)** and **Lemma 1**; t must be a divisor of p:
    $p = t \cdot u,\ u$ a positive integer.   **(45)**

From (45) and (44) we further get,
    $X \cdot Y = u \cdot r^2$   **(46)**
In virtue of $\gcd(X, Y) = 1$, equation (46) and **Proposition 1**

Imply that,
$$\begin{cases} X = u_1 \cdot r_1^2,\ Y = u_2 \cdot r_2^2,\ \text{and} \\ r = r_1 \cdot r_2.\ \text{Where } r_1, r_2, u_1 u_2 \\ \text{are positive integers such that} \\ \gcd(r_1, r_2) = 1 = \gcd(u_1, u_2), \\ \text{and } u = u_1 u_2 \end{cases} \quad (47)$$

Going back to **(42)** we see that (by **(47)**)
    $x = d \cdot X = d \cdot u_1 \cdot r_1^2,\ y = d \cdot Y = d \cdot u_2 \cdot r_2^2.$

And also by **(43)**, $z = t \cdot p = t \cdot (t \cdot u) = t^2 \cdot u = t^2 \cdot u_1 \cdot u_2$, by **(45)**
since $u = u_1 u_2$ by **(47)**

Finally $w = W \cdot d = r \cdot p \cdot d = r \cdot t \cdot u \cdot d = r_1 \cdot r_2 \cdot t \cdot u_1 \cdot u_2 \cdot d$
        by(43)            by(45)    $= d \cdot r_1 \cdot r_2 \cdot t \cdot u_1 \cdot u_2$

Also note that $\gcd(t, r_1 r_2) = 1$, since by (43), $\gcd(t, r) = 1$

And $r_1 r_2 = r$ ∎

## 11. Theorem 8 and its proof

We state Theorem 8, the last result of this work.

**Theorem 8**

*Consider the 5-variable (variables x,y,z,v,w) Diophantine system of two equations,*

$$\begin{cases} xy = v^2 \\ yz = w^2 \end{cases}$$

*Then, all postitive integer solutions of this system can be parametrically described in terms of 8-parameter formulas (parameters c,h,i,j,e,f,r,t),*

$$x = c \cdot h \cdot e^2 \cdot j^2 \cdot r^2, \; y = c \cdot h^3 \cdot e^2 \cdot j^4 \cdot i^2 \cdot f^2,$$
$$z = c \cdot h \cdot i^2 \cdot t^2, \; v = c \cdot i \cdot f \cdot h^2 \cdot e^2 \cdot j^3,$$
and $w = c \cdot e \cdot f \cdot t \cdot h^2 \cdot i^2 \cdot j^2$

*With the positive integer parameters $c, h, i, j, e, f, r, t$;*
*satisfying the conditions,*

$\gcd(i, j) = 1 = \gcd(e, f)$

**Proof**

A straightforward calculation shows that if five positive integers x,y,z,v, and w; satisfy the stated parametric formulas. Then the quintuple (x,y,z,v,w) is a solution to the above system. Conversely, below we prove that if (x,y,z,v,w) is a solution. Then the listed parametric formulas must be satisfied. We start with

$$\begin{cases} xy = v^2 \\ yz = w^2 \end{cases} \quad \textbf{(48)}$$

By Theorem 2 applied to each of the two equations in **(48)**; we must have,

$$\begin{cases} x = a \cdot r^2, \; y = a \cdot R^2, \; v = a \cdot r \cdot R; \text{ and} \\ z = b \cdot t^2, \; y = b \cdot T^2, \; w = b \cdot T \cdot t; \\ \text{where } a, r, R, b, T, t, \text{ are positive} \\ \text{integers with } \gcd(r, R) = 1 = \gcd(T, t) \end{cases} \quad \textbf{(49)}$$

Further, from (49) we must have,

$$a \cdot R^2 = b \cdot T^2 \quad (50)$$

Let $c = \gcd(a,b)$. So that,

$$\begin{cases} a = c \cdot a_1, \ b = c \cdot b_1 \\ \gcd(a_1, b_1) = 1 \end{cases} \quad (51)$$

From (50) and (51) we get,

$$a_1 \cdot R^2 = b_1 \cdot T^2 \quad (52)$$

Since $a_1$ and $b_1$ are relatively prime. Equation **(52)** and **Lemma2** imply that $a_1$ must be a divisor of $T^2$;

$$T^2 = a_1 \cdot k, \text{ for some positive integer } k \quad (53)$$

By (53) and Theorem 2 it follows that,

$$\begin{cases} a_1 = d \cdot e^2, k = d \cdot f^2; \text{ where} \\ d, e, \text{ and } f \text{ are positive integers} \\ \text{such that } \gcd(e, f) = 1 \text{ and } T = d \cdot e \cdot f \end{cases} \quad (54)$$

Combining (52) with (54) yields,

$$d \cdot e^2 \cdot R^2 = b_1 \cdot d^2 \cdot e^2 \cdot f^2;$$

$$R^2 = b_1 d f^2 \quad (55)$$

By Equation (55), $f^2$ is a divisor of $R$. Thus, by Lemma 2; f must be a divisor of R:

$$R = f \cdot g, \text{ for some positive integer } g \quad (56)$$

Thus from **(55)** and **(56)** we obtain,

$$g^2 = b_1 d \quad (57)$$

Equation (57) together with Theorem 2 imply that,

$$\begin{cases} b_1 = h \cdot i^2, \ d = h \cdot j^2, \text{ and} \\ g = h \cdot i \cdot j. \text{ where } h, i, j; \text{ are} \\ \text{positive integers with } \gcd(i, j) = 1 \end{cases} \quad (58)$$

We know make our way back. First from **(58)** and **(57)**
we have $R = f \cdot g = f \cdot h \cdot i \cdot j$

Next, from (54), $T = d \cdot e \cdot f = h \cdot j^2 \cdot e \cdot f$    by **(58)**

And $a_1 = d \cdot e^2 = h \cdot j^2 \cdot e^2$

Check to see that equation **(52)** is indeed satisfied:

$a_1 \cdot R^2 = (h \cdot j^2 \cdot e^2)(f \cdot h \cdot i \cdot j)^2 = h^3 \cdot j^4 \cdot e^2 \cdot i^2 \cdot f^2$

And $b_1 T^2 = (h \cdot i^2)(h \cdot j^2 \cdot e \cdot f)^2 = h^3 \cdot j^4 \cdot i^2 \cdot e^2 \cdot f^2$

Therefore from **(51)** and **(49)** we have,

$\quad x = ar^2 = a_1 \cdot c \cdot r^2 = c \cdot h \cdot j^2 \cdot e^2 \cdot r^2$

And $y = a \cdot R^2 = a_1 \cdot c \cdot R^2 = h \cdot j^2 \cdot e^2 \cdot c \cdot (f \cdot h \cdot i \cdot j)^2$;

$\quad y = c \cdot h^3 \cdot e^2 \cdot j^4 \cdot i^2 \cdot f^2$

Next, $z = b \cdot t^2 = c \cdot b_1 \cdot t^2 = c \cdot h \cdot i^2 \cdot t^2$

And
$\quad v = a \cdot r \cdot R = c \cdot a_1 \cdot f \cdot h \cdot i \cdot j \cdot r;$
$\quad v = c \cdot (h \cdot j^2 \cdot e^2) h \cdot f \cdot i \cdot j \cdot r = c \cdot i \cdot f \cdot h^2 \cdot e^2 \cdot j^3$

Finally,
$\quad w = b \cdot T \cdot t = c \cdot b_1 \cdot T \cdot t = c \cdot (h \cdot i^2)(h \cdot j^2 \cdot e \cdot f) \cdot t;$
$\quad w = c \cdot e \cdot f \cdot t \cdot h^2 \cdot i^2 \cdot j^2$

And with (from **(54)** and **(58)**), $\gcd(i, j) = 1 = \gcd(e, f)$. ∎

## 12. Geometric interpretations

An obvious geometric interpretation of the Diophantine equation, $xy = z^n$; with n a fixed positive, $n \geq 2$, is the following: The set of all positive integer solutions to such an equation in effect describes the entire family of integer-sided rectangles whose area is an nth integer power. On the other hand, the solutions to the equation, $xyz = w^2$. Can be thought to describe the family of all rectangular parallelepipeds with integer side lengths, and whose volume is a perfect square. Finally, the geometric interpretation to the Diophantine system in Theorem 8; refers to the family of rectangular parallelepipeds each containing two pairs (each pair consisting of two opposite and congruent rectangles) of lateral faces. With each face (rectangle) having area which is an integer square. (A simpler term, instead of parallelepipeds; might be, 3-D boxes)